\documentclass{amsart} 
\usepackage{amssymb}

\title{How to release Frege's system from Russell's antinomy} 
\author{Paola Cattabriga}
\address{University of Bologna --  Italy}
\email{co14099@iperbole.bologna.it}
\begin{document}
\begin{abstract}
The conditions for proper definitions in
mathematics are given, in terms of the  theory of definition, on the basis 
of the criterions of eliminability and
non-creativity. As a definition, Russell's antinomy is a violation of the criterion
of eliminability (Behmann, 1931, \cite{behmann};  Bochvar, 1943, \cite{bochvar}).  
Following the path
of the criterion  of
non-creativity, this paper develops a new analysis of Comprehension schema
 and, as a consequence, proof that Russell's
antinomy argumentation, despite the words of Frege himself, does not hold 
in Grundgesetze der Arithmetik.  According to Basic Law 
(\ref{blIII}), the class of 
classes not belonging to themselves is a class 
defined by a function which can not take as argument its own course of 
value, 
\begin{equation}\notag 
 \put(0,2.3){\line(1,0){8}}\put(4,0){\line(0,1){2.3}}\quad \bigg ( 
{\boldsymbol\forall}
 = \grave{\varepsilon}\Big ( 
\put(-3.4,3.5){\line(1,0){3.4}}\put(0.1,0.9){\line(0,4){2.3}}\overset{\frak{g}}{\put(-0.2,3.5){\line(1,0){3.5}}\smile\put(-3.4,3.5){\line(1,0){3.5}}}
\big ( \, \grave{\varepsilon}(\overline{\quad\; }\frak{g}(\varepsilon))= 
\varepsilon \,  \rightarrow \frak{g}(\varepsilon)
\big ) 
\Big ) \bigg ),
\end{equation} 
in other words, the class of classes not belonging to themselves is a class 
whose classes are not identical to the class itself \cite{lc06}.
\end{abstract}
\maketitle
\thispagestyle{empty}

\section*{Introduction}
Recently, a considerable attention has been devoted to the discovery of 
new ways for repairing Frege's system from the damage produced by 
Russell's antinomy. 
All these 
attempts share a common aptitude in finding new restrictive forms or 
applications of the comprehension principle  
\cite{burgess,heck,parsons,sh,wehmeier}, either  offer a  
revised interpretation of Basic Law (V) \cite{boolos}.

The present note proposes a  completely new way to restore 
Frege's system, neither restrictive nor revising,  rising from the logical
inquiry in the theory of definition and the author's 
investigations in self-reference procedures \cite{cA,cC}.
Partially in agreement with their constructive 
attitude, the result we are concerned with draws its first inspiration 
from  
Jules-Henri Poincar\'e's writings about the significance and riskiness of 
impredicative 
definitions in mathematics.
\begin{quote}
{ \sl 
The fact is that such a procedure is not   applicable. Why? 
Because their definitions
are not predicative and  contain within such a vicious circle  I already 
mentioned  
above; not predicative definitions can not be \emph{substituted to 
defined terms}. In this condition,
logistics is no longer sterile: it generates contradictions. }
(1902, \cite{poincare}, 211, our translation) 
\end{quote}
With these words Poincar\'e was directly referring to Cantor's famous
diagonalization argument for reaching more-than-denumerable sets, 
criticizing the underlying belief of the  existence of the
actual infinite.  Russell's antinomy and Cantor's diagonalization are both
self-referring procedures, and by admission of Russell himself his 
discovery is due to his reflections about Cantor's argument.
\begin{quote}
{ \sl I was led to the contradiction in the following way. As you, of course 
know, Cantor proved that there is no greatest number. [$\dots$] Now there 
are concepts 
whose extension comprises everything; these should therefore have the 
greatest number. I tried to set up a one-one relation between all objects 
and all classes; when I applied Cantor's proof with my special relation, I 
found that the class $Cls\cap x\backepsilon (x\sim\epsilon x)$ 
was left over, even though all classes had already been 
enumerated. I have already been thinking about this contradiction for a 
year; 
I believe the only solution is that function and argument must be able to 
vary independently.} (\cite{correspondence}, 215, XXXV/3 Russell 
to Frege 24.6.1902)
\end{quote}
Cantor's diagonalization and Russell's antinomy are logically  
connected, and for the historical perspective the latter rises from 
the former. 
We shall expose how the logical inquiry regarding the mathematical 
definitions and the rules for applying substitutions, as intuited by 
Poincar\'e,  leads to the restore of Frege's system. 
Nowadays Russell's antinomy is widespread in many several forms 
throughout logic and mathematics,  we shall apply 
our argument to three different theories.
First, in section \ref{zf}, to Zermelo-Fraenkel set theory (\emph{ZF}) 
\cite{foundations}, and to  
a first order set theory close to Frege's system 
(\cite{hatcer} 76-90),
and finally in section \ref{fgga}, directly to Frege's Grundgesetze der 
Arithmetik  \cite{gga}. 
For this last section the redear is required  to be familiar  with the  original 
notation of  Grundgesetze.
This paper is the  revised and extended version of \cite{cB}.
\section{ZERMELO-FRAENKEL SET THEORY}\label{zf} 
In Cantor's \emph{naive} set theory a set was 
thought to be an aggregate, i.e  a ``collection of elements into a 
whole of our intuition".  This conception appears  accordingly 
symbolized by the  principle of comprehension, which can be regarded as a 
schema
originally conceived for defining sets.
For this special function, we are legitimate in 
considering the comprehension principle under the point of view of the 
theory of 
definition\footnote{
We are referring here to the field of logical researches of ancient 
tradition, which in the last century
has seen  noteworthy contributions from A. Padoa, 
J.C.C. Mckinsey, A. Tarsky,  E. W. Beth, A. Robinson, W. Craig.  
It is also worth  mentioning 
how  significative the theory of definition was throughout Frege's whole
foundation  of mathematics  \cite{beg,gla,gga}, even if it was unfortunately 
neglected, as we will see, in the 
famous rushed appendix (\cite{gga}, II 253-265). The principle of 
Comprehension  results to be  formalized for the first time in  \emph{Grundgesetze},  
$$ \vdash f(a) = a\smallfrown  \grave{\varepsilon}f(\varepsilon),$$ but its role in defining 
\emph{new objects}  seems not to be sufficiently focussed by Frege
(\cite{gga}, I 73, 75, 117). All that appears to be comprehensible for that 
time, as the modern theory of definition was at its beginning  exactly with 
Frege and Peano.
Remained in the twilight for 
years, the first notable result investigating  the meaning of the introduction of 
\emph{new symbols} for the interplay between mathematical definition 
and logical deduction,  could be viewed as an outcome of the questions raised in 
Frege--Peano correspondence \cite{padoa}. 
}. 
 On that basis the rules for  definitions,  established by the criterions 
of eliminability 
and non-creativity,  state the conditions for proper equivalence, giving 
some basic 
restrictions to prevent superimpositions and circularity (\cite{suppes}, 
151-173).
The rule for defining a \emph{new operation symbol} (or a new individual 
constant, i.e. an operation symbol of rank zero) is as follows:
\begin{itemize}
\item[$r_{1}$]An equivalence $E$ introducing a new 
n-place operation symbol $O$ is a 
proper definition in a theory if, and only if, $E$ is of the form 
$$ O(x_{1}\dots x_{n}) = y \leftrightarrow \psi,$$ and the following restrictions are 
satisfied: \emph{(i)} $x_{1}\dots x_{n},$ $ y$ are distinct variables, 
\emph{(ii)} 
$\psi$ has no free variables other than $x_{1}\dots x_{n}, y$, 
\emph{(iii)} $\psi$ is a formula in which the only non-logical constants 
are 
primitive symbols and previously defined symbols of the theory, 
and \emph{(iv)} the formula $\exists ! y \psi$
is derivable from the axioms and preceding definitions of the theory.
\end{itemize}
With condition \emph{(iv)} we are required to have a preceding theorem 
which guarantees that the  operation is uniquely defined.
If the 
restriction on the uniqueness is dropped then a contradiction can be 
derived  (\cite{suppes}, 159). 
The self-referring characteristic of 
Russell's antinomy  
will turn  out to be a procedure which forces a set to be twice defined, 
we are 
dealing with a set 
whose elements belong and do not belong to the set itself, dropping thus 
the above uniqueness condition.

As it is well-known, in \emph{ZF} when defining a set, 
the uniqueness for the comprehension axiom schema  is ensured 
by the axiom of extensionality. Let us introduce to them both,
\begin{align}\tag{$C$}\label{C}
\text{\emph{Comprehension }} &\; \forall z_{1}...\forall z_{n}\exists y\forall 
x(x \in y \leftrightarrow \varphi(x)),\\
\intertext{where $ \varphi (x)$ is any formula in \emph{ZF},  $z_{1},...,z_{n}$ are the free variables of $\varphi (x)$ other than $x$, 
and $y$ is not a free variable of $\varphi(x)$,}
\tag{$E_{1}¥$}\label{E1}
\text{\emph{Extensionality }} &\; \forall x \forall y[ \forall z(z \in x \leftrightarrow z \in y) 
\leftrightarrow x = y].\end{align}
\emph{Extensionality}, without any additional axioms,  implies that for 
every condition $\varphi(x)$ on $x$ (in \emph{Comprehension}) there exists ``at 
most" one 
set $y$ which contains \emph{exactly} those elements $x$ which fulfil the 
condition 
$\varphi(x)$. 
In other words, if  there is a set $y$ such that $\forall 
x(x \in y \leftrightarrow \varphi(x))$, $y$ is unique. It can be shown as 
follows. If $y'$ is also such, i.e. $\forall 
x(x \in y' \leftrightarrow \varphi (x))$, then we have, obviously $\forall 
x(x \in y' \leftrightarrow x \in y)$, and then by \emph{Extensionality}, $y' = y$. 
(\cite{foundations} 31, \cite{levy} 7, \cite{bourb} 67). 

Usually Russell's antinomy argumentation is presented as 
follows.
\begin{enumerate}
\item[\it RA ]
\emph{There exists no set  which 
contains exactly those elements which do not contain 
themselves, in symbols  $\neg\exists y\forall x(x\in y\leftrightarrow x\notin 
x) $.}
\item[]Proof. By contradiction. Assume that $y$ is a set such that for every 
element 
$x$, $ x \in y$ if and only if $ x \notin x$. For $x=y$, we have $y \in y$
 if and only 
if $y \notin y$. Since, obviously, $y \in y$ or $y \notin y$, and as we 
saw, 
each of $y \in y$ and $y \notin y$ implies the other statement, we have 
both $y 
\in y$ and $y \notin y$, which is a contradiction (\cite{foundations}, 31).
\end{enumerate}
This proof holds for sure as a first order theorem, as indeed in first 
order logic $$\vdash \neg\exists y\forall x(x\in y\leftrightarrow x\notin 
x),$$
and it can affect 
comprehension principle when it is regarded individually or in itself
$$(  \text{\ref{C}} ) \rightarrow \exists y\forall x(x\in y\leftrightarrow x\notin 
x) \vdash \, ,$$
 but 
it does not hold in 
any system which applies both \emph{Comprehension }and \emph{Extensionality}
$$(  \text{\ref{C}} ) \cup (  \text{\ref{E1}} ) \vdash \neg \forall x(x\in y\leftrightarrow x\notin 
x).$$  
Indeed we can show that Russell's antinomy does not affect a first 
order set theory, since by \emph{Extensionality} the 
above proof \it RA \rm can not be accomplished. If as an example of   
\emph{Comprehension} we define
$$x \in y \text{ if and only if } x \notin x$$
then by \emph{Extensionality}   we obtain always (detailed proof 
 in the following of this section and in section \ref{fgga})
\begin{equation}\tag{$*$}\label{star}
(x \in y \leftrightarrow x \notin x) \rightarrow \neg(x = y),
\end{equation}
so it is never 
the case that 
$x = y.$ 
This is an applicative case of the uniqueness 
condition.
Indeed a special case because it involves negation and therefore 
complementation.
For it, the above inference  \it RA \rm
can not be accomplished since $``x = y"$ can not be assumed and $``y \in y 
\text{ if and only if }  y \notin y"$ is not derived.
In other terms, for a simple first order rule from $$\forall x(x \in a \leftrightarrow 
x \notin x)$$ we can yield $$(a \in a \leftrightarrow a \notin a),$$  but this rule is not 
applicable during a reasoning  
of set theory where, by \emph{Extensionality}, the formula 
\begin{equation}\tag{$**$}\label{star1}
(x \in a \leftrightarrow x \notin x) \rightarrow \neg(x = a)
\end{equation}
can always be derived.

Two sets are 
equal when  and only  when they have the same elements, and the equality between sets
must fulfill all the requirements of the relation of 
identity, 
namely, reflexivity,  every set has the same elements as itself, and 
substitutivity of identity, $ x = y \rightarrow 
(P(x) \leftrightarrow P(y))$  (for $x$ and $y$  any set and $P$ 
 any property). Actually the substitutivity of identity is a version of 
the Principle of the Indiscernibility of 
Identicals, one of the consequences of Frege's Basic Law 
(\ref{blIII})   (see  section  \ref{fgga}). 
In set theory it doesn't follow directly from the relation of 
coexstensiveness, and, whether  axiom or   
hypothesis, a principle must be stated
that two coextensive sets 
 share all their properties and are therefore
equal. 
The principle of extensionality says in effect that two 
sets are identical if and only if they have the same elements, in 
symbols (\ref{E1}). 
Let us  recall also a further first order schema 
as follows,  
\begin{align}\tag{$E_{2}¥$}\label{E2}
\text{\emph{Extensionality }} &\; \forall x \forall y[  x = y \rightarrow  
(\phi (x,x) \leftrightarrow \phi (x,y)) ],
\end{align}
where $\phi (x,y)$ is obtained from $\phi (x,x)$ by replacing $y$ for zero, one or more 
occurrences of $x$ in the wff $\phi (x,x)$, and $y$ is free for $x$ in all 
occurrences of $x$ which it replaces (\cite{hatcer} 79). Under the logical point of view (\ref{E2}) is,
as respect to (\ref{E1}), the integral 
 schema, quite close to
 Frege's Basic Law (\ref{blIII}). 
(\ref{E2}) yields
$$\forall x \forall y[  x = y\rightarrow \forall z(z \in x \leftrightarrow 
z \in y) ]$$
 (\cite{hatcer} 139), which  holds usually in \emph{ZF} together 
  with the further statement
$$\forall x \forall y[  \forall z(z \in x \leftrightarrow z \in y) \rightarrow 
x = y],$$    
either as axiom (\cite{foundations} 27), or theorem (\cite{hatcer} 141).
Any derivation of (\ref{star}) in \emph{ZF} requires necessarily both 
statements, so that we adopted directly (\ref{E1}) (\cite{lavrov} 109).
As in Frege's Grundgesetze, see section (\ref{fgga}), (\ref{star}) follows 
easily from (\ref{E2}), whereas  the derivation of (\ref{star})  from (\ref{E1}) 
is bound to the methodological concerns of 
\emph{ZF}. We shall expose  both the cases. 

Let us begin with the first,
\begin{equation}\tag{\sl I}\label{E2-1}
x = y \rightarrow  (x \in x \leftrightarrow x \in y )
\end{equation}
is an instance of (\ref{E2}) where $\phi (x,x)$ is $ x \in x$ and 
$\phi (x,y)$ is $ x \in y$, so that $ x \in 
y $ is obtained by substituting $x$ with $y$ in $x \in x$. We notice 
that the restriction on $y$ in (\ref{E2}) is 
fulfilled so that (\ref{E2-1}) follows logically from (\ref{E2}).
(\ref{E2-1}) yields by \small $(p \rightarrow q) \leftrightarrow (\neg q \rightarrow \neg 
p)$\normalsize  and \small $(p \leftrightarrow q) \leftrightarrow (q 
\leftrightarrow  p)$\normalsize 
\begin{equation}\tag{\sl II}\label{E2-2}
\neg (x \in y \leftrightarrow x \in x ) \rightarrow \neg ( x = y)
\end{equation}
and then,\small $\neg ( p \leftrightarrow q) \leftrightarrow ( p \leftrightarrow \neg  q)$,\normalsize 
\begin{equation}\tag{\sl III}\label{E2-3}
(x \in y \leftrightarrow x \notin x ) \rightarrow  \neg ( x = y),
\end{equation}
i.e. (\ref{star}) QED.

In \emph{ZF} (\ref{star})
can be  derived easily from (\ref{E1}) too, with the exception of  the 
methodological limitations about illegitimate self-referring.
\begin{quote}
{ \sl
In most cases when we use variables as metamathematical variables for 
variables we shall assume, tacitly, that different variables stand for 
different variables. E.g. the set of all statements 
$\forall z ( z \in x  \leftrightarrow   z \in y )$
is also assumed to contain the statement $\forall u ( u \in w  \leftrightarrow   u 
\in v )$, but not to contain the statement $\forall x ( x \in x  \leftrightarrow   
x \in y )$. In some other cases we do not insist that different 
variables stand for different variables. 
} (\cite{foundations} 21)
\end{quote}
As well-known set theorists use methodically such tacit assumptions  in 
\emph{ZF} with the purpose of preventing 
self-contradictory procedures from arising. 
But our aim is actually just to show that Russell's antinomy does not 
follow when  
 \emph{Extensionality} is taken suitably into account.
 We shall place therefore the above  
tacit assumptions in stand by, and we shall keep  in enquiring instead into  the restrictions established 
 by the theory of definition.

 Furthermore, let us add, that
if we  allow expressions of totalities and self-reference in 
\emph{Comprehension} schema there is no reason to prevent the same 
expressions within \emph{Extensionality} 
schema.

 Accordingly, in \emph{ZF} 
$$  (x \in y \leftrightarrow x \in x) \leftrightarrow x = y,$$ is an 
instance of   (\ref{E1}),    hence
$$ \neg (x \in y \leftrightarrow x \in x) \rightarrow \neg( x = y), $$ 
and finally
$$ (x \in y \leftrightarrow x \notin x) \rightarrow \neg( x = y),$$
i.e. (\ref{star}), QED.
Whenever (\ref{C}) yields $\forall x (x \in a \leftrightarrow x \notin 
x)$ by (\ref{star})
we have     (\ref{star1}) too.

We can now state that Russell's antinomy argumentation does not hold
in a first order system based on \emph{Comprehension} and 
\emph{Extensionality}, because it violates
the restriction on the uniqueness established by extensionality.
In defining a set containing exactly those elements which do not 
contain themselves we define a set which members can not be considered to 
be
identical to the set itself, either we should have a set defined twice, as 
a set which contains and does not contain itself. This would be to 
define a set and its complement as being exactly the same identical 
set, and in detail to violate the criterion of non-creativity and its 
consequent criterion of relative consistency (\cite{suppes}, 155).
From the logical point of view  Russell's antinomy 
turns out to be just a partial argument. Although set theory disposes of the 
symbols to mention, or to express,  the existence of  sets containing 
exactly those elements which do not 
contain themselves, yet, when the argument reaches its own whole 
representation,
no contradiction can be derived. 
All that  turns out to be amazing especially in Frege's Grundgesetze, 
since the subcomponent 
``$ 
\put(0,2.3){\line(1,0){8}}\put(4,0){\line(0,1){2.3}}\quad 
({\boldsymbol\forall} = {\boldsymbol\forall} )$",
which prevent the contradiction in Frege's \emph{way out}  (\cite{gga}, 
II 262-263),
can be  derived  simply by Basic Law (\ref{blIII}), namely 
without Frege's restrictions. The so called ``Frege's \emph{way out}" is in 
the relative literature the shortening for the original Frege's attempt to 
repair his system, by means of
(V'b) and (V'c). As well-known  the amended axiom blocks the derivation of 
Russell's antinomy, but it is not formally 
satisfactory (\cite{sluga} 165-170).  On the contrary the solution we 
propose
here  seems to be formally faultless.
\section{FREGE SYSTEM}\label{fgga} 
We are now ready to 
reconsider the  appendix of 
Grundgesetze der Arithmetik (\cite{gga}, II 253-265) in accordance with what was previously stated. We shall prove that 
despite the words of Frege himself Russell's antinomy can not be derived.

In the following, Frege's symbolism is fully adopted, except for the 
modern symbol
of implication.
First let us look at the deduction involved with function 
$\zeta\smallfrown\xi$.
\begin{quote}
{ \sl
 
... Now let us see how the matter turns out if we make use of our sign 
``$\smallfrown$". 
Here ``$ \grave{\varepsilon}( 
\put(0,2.3){\line(1,0){8}}\put(4,0){\line(0,1){2.3}}\quad\varepsilon\smallfrown \varepsilon )$" 
will occur in place of  ``${ \boldsymbol\forall }$". 
If in our proposition 
\begin{equation}\tag{$82$}
F(a \smallfrown \grave{\varepsilon}f(\varepsilon)) \; 
\rightarrow \;  F(f(a))
\end{equation}
we take 
\begin{align}
\tag*{$i)$}
`` \put(0,2.3){\line(1,0){8}}\put(4,0){\line(0,1){2.3}}\quad\xi\smallfrown 
\xi"
\text{ for } ``f(\xi)", &
\\  
\tag*{$ii)$}\label{sub2}
 ``\overline{\quad\; } \xi" \text{ for }  ``F(\xi)", &
\\ 
\tag*{$iii)$}  \text{and } ``\grave{\varepsilon}( 
\put(0,2.3){\line(1,0){8}}\put(4,0){\line(0,1){2.3}}\quad\varepsilon
\smallfrown \varepsilon)"  \text{ for }  ``a", &
\end{align}
then we obtain
\begin{multline}\tag{$\theta$}\label{teta}
\grave{\varepsilon}( 
\put(0,2.3){\line(1,0){8}}\put(4,0){\line(0,1){2.3}}\quad\varepsilon\smallfrown 
\varepsilon)  
\smallfrown \grave{\varepsilon}( 
\put(0,2.3){\line(1,0){8}}\put(4,0){\line(0,1){2.3}}\quad\varepsilon\smallfrown 
\varepsilon) 
 \rightarrow \\
\,
 \text{ } \put(0,2.3){\line(1,0){8}}\put(4,0){\line(0,1){2.3}}\quad
\grave{\varepsilon}( 
\put(0,2.3){\line(1,0){8}}\put(4,0){\line(0,1){2.3}}\quad\varepsilon\smallfrown 
\varepsilon) 
\smallfrown \grave{\varepsilon}( 
\put(0,2.3){\line(1,0){8}}\put(4,0){\line(0,1){2.3}}\quad\varepsilon\smallfrown 
\varepsilon), 
\end{multline}
from which by (Ig) there follows
\begin{equation}\tag{$\iota$}\label{iota}
\put(0,2.3){\line(1,0){8}}\put(4,0){\line(0,1){2.3}}\quad
\grave{\varepsilon}( 
\put(0,2.3){\line(1,0){8}}\put(4,0){\line(0,1){2.3}}\quad\varepsilon\smallfrown 
\varepsilon)  
\smallfrown \grave{\varepsilon}( 
\put(0,2.3){\line(1,0){8}}\put(4,0){\line(0,1){2.3}}\quad\varepsilon\smallfrown 
\varepsilon).
\end{equation}
Making the same substitutions in proposition 
\begin{equation}\tag{$77$}
F(f(a))\;\rightarrow\; F(a \smallfrown 
\grave{\varepsilon}f(\varepsilon))
\end{equation}
we obtain
\begin{multline}\tag{$\kappa$}\label{kappa}
\put(0,2.3){\line(1,0){8}}\put(4,0){\line(0,1){2.3}}\quad
\grave{\varepsilon}( 
\put(0,2.3){\line(1,0){8}}\put(4,0){\line(0,1){2.3}}\quad\varepsilon\smallfrown 
\varepsilon)  
\smallfrown \grave{\varepsilon}( 
\put(0,2.3){\line(1,0){8}}\put(4,0){\line(0,1){2.3}}\quad\varepsilon\smallfrown 
\varepsilon) 
 \rightarrow  \\
\,
 \grave{\varepsilon}( 
\put(0,2.3){\line(1,0){8}}\put(4,0){\line(0,1){2.3}}\quad\varepsilon\smallfrown 
\varepsilon)  
\smallfrown \grave{\varepsilon}( 
\put(0,2.3){\line(1,0){8}}\put(4,0){\line(0,1){2.3}}\quad\varepsilon\smallfrown 
\varepsilon),
\end{multline}
from which together with ($\iota$) there follows
\begin{equation}\tag{$\lambda$}\label{lambda}
 {\overline{\quad\; }\,}
\grave{\varepsilon}( 
\put(0,2.3){\line(1,0){8}}\put(4,0){\line(0,1){2.3}}\quad\varepsilon\smallfrown 
\varepsilon)  
\smallfrown \grave{\varepsilon}( 
\put(0,2.3){\line(1,0){8}}\put(4,0){\line(0,1){2.3}}\quad\varepsilon\smallfrown 
\varepsilon), 
\end{equation}
which contradicts ($\iota$). Therefore at least one of the two 
(77) and (82) must be false, and therefore proposition 
\begin{equation}
 \vdash f(a) = a\smallfrown\grave{\varepsilon}f(\varepsilon),
\end{equation} 
also, from 
which they both follow. A look at the derivation of  (1) in \S 55 of 
our first 
volume shows that there too use is made of (Vb). Thus suspicion is 
directed at this proposition here as well. } (\cite{gga}, II 257) 
\end{quote}
Within Frege's system the proposition $$ \vdash a=b \, \rightarrow \, (g(a) = 
g(b))$$ (IIIh, \cite{gga}, I 66-67) can be achieved by basic law 
\begin{equation}\tag{III}\label{blIII}
\vdash f(a = b) \rightarrow f( 
\overset{\frak{g}}{\put(-0.2,3.5){\line(1,0){3.5}}\smile
\put(-3.4,3.5){\line(1,0){3.5}}}  ( \frak{g}(b) \rightarrow 
\frak{g}(a))),
\end{equation}
as well as the proposition 
$(\put(0,2.3){\line(1,0){8}}\put(4,0){\line(0,1){2.3}}\quad g(a) = 
g(b) )\; \rightarrow \;  
\put(0,2.3){\line(1,0){8}}\put(4,0){\line(0,1){2.3}}\quad (a = b).$

In detail from 
$$\,  \put(0,2.3){\line(1,0){8}}\put(4,0){\line(0,1){2.3}}\quad g(a) \, \rightarrow \, 
(g(b) \,\rightarrow 
\,\put(0,2.3){\line(1,0){8}}\put(4,0){\line(0,1){2.3}}\quad(a=b))$$
(IIIb, ibid. 66-67), and
$g(\xi)=( \,\put(0,2.3){\line(1,0){8}}\put(4,0){\line(0,1){2.3}}\quad 
g(a)=g(\xi)) $ we obtain 

$$
\,\put(0,2.3){\line(1,0){8}}\put(4,0){\line(0,1){2.3}}\text{ }
\,\put(0,2.3){\line(1,0){8}}\put(4,0){\line(0,1){2.3}}
\quad g(a) =g(a)
\,\rightarrow\, ((\put(0,2.3){\line(1,0){8}}\put(4,0){\line(0,1){2.3}}\quad 
g(a)=g(b)) \,\rightarrow\, 
\put(0,2.3){\line(1,0){8}}\put(4,0){\line(0,1){2.3}}\quad(a=b)),
$$
and
$$ g(a) =g(a)
\,\rightarrow\, ((\put(0,2.3){\line(1,0){8}}\put(4,0){\line(0,1){2.3}}\quad 
g(a)=g(b)) \,\rightarrow\, 
\put(0,2.3){\line(1,0){8}}\put(4,0){\line(0,1){2.3}}\quad(a=b)),$$
 so that, by  $\vdash a = a $ (IIIe, ibid. 66-67),
\begin{equation}\label{ext}
(\put(0,2.3){\line(1,0){8}}\put(4,0){\line(0,1){2.3}}\quad g(a) = 
g(b) )\; \rightarrow \;  
\put(0,2.3){\line(1,0){8}}\put(4,0){\line(0,1){2.3}}\quad (a = b).
\end{equation}

Let us now reconsider in accordance  the above Frege's substitutions in (82) 
and (77). Firstly, 
by the substitution $i)$ in law (1)  we obtain
\begin{equation}\label{c1}
 \put(0,2.3){\line(1,0){8}}\put(4,0){\line(0,1){2.3}}\quad a \smallfrown a
\; = \;
 a \smallfrown \grave{\varepsilon}( 
\put(0,2.3){\line(1,0){8}}\put(4,0){\line(0,1){2.3}}\quad
\varepsilon \smallfrown \varepsilon ).
\end{equation}
Then  replacing 
\begin{align}
 ``a  \smallfrown \xi " & \text{ for }  
``g(\xi)", \tag*{$iv)$}\label{subIIIb1} \\ 
\text{and } 
``\grave{\varepsilon}( 
\put(0,2.3){\line(1,0){8}}\put(4,0){\line(0,1){2.3}}\quad
\varepsilon \smallfrown \varepsilon )" & \text{ for }  ``b", 
\tag*{$v)$}\label{subIIIb2}
\end{align}
in ($2$) we obtain 
\begin{equation}\label{anti1}
(\;\put(0,2.3){\line(1,0){8}}\put(4,0){\line(0,1){2.3}}\quad a \smallfrown a\; =
 a \smallfrown \grave{\varepsilon}( 
\put(0,2.3){\line(1,0){8}}\put(4,0){\line(0,1){2.3}}\quad
\varepsilon \smallfrown \varepsilon )  ) 
\; \rightarrow  \;
\put(0,2.3){\line(1,0){8}}\put(4,0){\line(0,1){2.3}}\quad ( a =
\grave{\varepsilon}( 
\put(0,2.3){\line(1,0){8}}\put(4,0){\line(0,1){2.3}}\quad
\varepsilon \smallfrown \varepsilon )).
\end{equation}
Therefore by ($3$) and ($4$)
\begin{equation}
\put(0,2.3){\line(1,0){8}}\put(4,0){\line(0,1){2.3}}\quad (a =
\grave{\varepsilon}( 
\put(0,2.3){\line(1,0){8}}\put(4,0){\line(0,1){2.3}}\quad
\varepsilon \smallfrown \varepsilon )).
\end{equation}
Now the replacement of
``$\grave{\varepsilon}( 
\put(0,2.3){\line(1,0){8}}\put(4,0){\line(0,1){2.3}}\quad
\varepsilon\smallfrown \varepsilon)$" 
for ``$a$" in  ($3$)  turns out to be invalid, and consequently
Frege's substitution $iii)$ in (82) and (77)
 can not be accomplished. Therefore ($\theta$),  
($\iota$), ($\kappa$),  ($\lambda$) and the contradiction can no 
longer be derived\footnote{Compare with Frege's amendment, $(1'$ (\cite{gga} 
264-265)}, QED.

Our proof shows that Frege's \emph{way out} of negating for any function the 
possibility to take as argument  its own course of value 
(\cite{gga}, II  262),  as in 
\begin{equation}\tag{V'b}\label{V'b}
\vdash  (\grave{\varepsilon}f (\varepsilon) = \grave{\alpha}g(\alpha)) 
\rightarrow  
( \put(0,2.3){\line(1,0){8}}\put(4,0){\line(0,1){2.3}}\quad(a =  
\grave{\varepsilon}f (\varepsilon))  
\rightarrow  f(a) = g(a))
\end{equation}
and 
\begin{equation}\tag{V'c}\label{V'c}
\vdash  (\grave{\varepsilon}f (\varepsilon) = \grave{\alpha}g(\alpha)) 
\rightarrow  
( \put(0,2.3){\line(1,0){8}}\put(4,0){\line(0,1){2.3}}\quad(a =  
\grave{\alpha}g(\alpha)) \rightarrow   f(a) = g(a)),
\end{equation}
turns out to be useless and too restrictive.
Frege's Basic Law (III) already
 rules out  the cases of those functions which would lead to contradiction 
taking as argument its own course of value, without need to exclude all 
the 
other cases of self-referring. Indeed there can be cases of functions 
which argument 
can be its own course of value without yielding any contradiction, so that
there is no reason for excluding them a priori.
To clarify all that, let us recall Frege's passage after the 
introduction of (V'b) and (V'c).
\begin{quote}
{ \sl

Let us now convince ourselves that the contradiction that arose earlier 
between the proposition ($\beta$) and ($\epsilon$) is now avoided. We 
proceed as we did in the derivation of ($\beta$), using (V'c) 
instead of (Vb). As before, let 
\begin{equation}\notag
 *) \text{ } 
``{\boldsymbol\forall}" \text{ abbreviate } 
``\grave{\varepsilon}\Big
 (
\put(-3.4,3.5){\line(1,0){3.4}}\put(0.1,0.9){\line(0,4){2.3}}
\overset{\frak{g}}{\put(-0.2,3.5){\line(1,0){3.5}}\smile\put(-3.4,3.5){\line(1,0){3.5}}}
\big ( \, \grave{\varepsilon}(\overline{\quad\; }\frak{g}(\varepsilon))= 
\varepsilon \,  \rightarrow \frak{g}(\varepsilon)
\big ) \Big )".
\end{equation}
By  (V'c) we have   
\begin{multline}\notag
{\boldsymbol\vdash} \grave{\varepsilon}(\overline{\quad\; 
}f(\varepsilon))
 =  
\grave{\varepsilon}\Big (  \put(-3.4,3.5){\line(1,0){3.4}}\put(0.1,0.9)
{\line(0,4){2.3}}\overset{\frak{g}}{\put(-0.2,3.5){\line(1,0){3.5}}\smile\put(-3.4,3.5){\line(1,0){3.5}}}
\big ( \grave{\varepsilon}(\overline{\quad\; }\frak{g}(\varepsilon))= 
\varepsilon   \rightarrow \frak{g}(\varepsilon)
\big ) 
\Big )
\rightarrow 
\\
\Bigg ( 
\put(0,2.3){\line(1,0){8}}\put(4,0){\line(0,1){2.3}}\quad \bigg 
({\boldsymbol\forall}
= 
\grave{\varepsilon}\Big (  \put(-3.4,3.5){\line(1,0){3.4}}\put(0.1,0.9)
{\line(0,4){2.3}}\overset{\frak{g}}{\put(-0.2,3.5){\line(1,0){3.5}}\smile\put(-3.4,3.5){\line(1,0){3.5}}}
\big ( \, \grave{\varepsilon}(\overline{\quad\; }\frak{g}(\varepsilon))= 
\varepsilon   \rightarrow \frak{g}(\varepsilon)
\big ) 
\Big ) \bigg ) 
\rightarrow  \\
 (\overline{\quad\; }f({\boldsymbol\forall})) = \text{ }
\put(-3.4,3.5){\line(1,0){3.4}}\put(0.1,0.9){\line(0,4){2.3}}\overset{\frak{g}}
{\put(-0.2,3.5){\line(1,0){3.5}}\smile\put(-3.4,3.5){\line(1,0){3.5}}}
\big ( \,
 \grave{\varepsilon}(\overline{\quad\; }\frak{g}(\varepsilon))= 
{\boldsymbol\forall}
   \rightarrow \frak{g}({\boldsymbol\forall})
\big ) 
\Bigg ).
\end{multline} 
Using our abbreviation, we obtain  
\begin{multline}\notag
 {\boldsymbol\vdash} \grave{\varepsilon}(\overline{\quad\; 
}f(\varepsilon))  = 
\grave{\varepsilon}\Big (  
\put(-3.4,3.5){\line(1,0){3.4}}\put(0.1,0.9){\line(0,4){2.3}}\overset{\frak{g}}
{\put(-0.2,3.5){\line(1,0){3.5}}\smile\put(-3.4,3.5){\line(1,0){3.5}}}
\big ( \,
 \grave{\varepsilon}(\overline{\quad\; }\frak{g}(\varepsilon))= 
\varepsilon \,  \rightarrow \frak{g}(\varepsilon)
\big ) 
\Big )
\rightarrow 
\\
\Bigg ( \put(0,2.3){\line(1,0){8}}\put(4,0){\line(0,1){2.3}} 
\quad({\boldsymbol\forall}
= 
{\boldsymbol\forall} )
\rightarrow \quad\quad\quad\quad\quad\quad\quad\quad\quad
\\
(\overline{\quad\; }f({\boldsymbol\forall}))
 = \text{  }
\put(-3.4,3.5){\line(1,0){3.4}}\put(0.1,0.9){\line(0,4){2.3}}\overset{\frak{g}}
{\put(-0.2,3.5){\line(1,0){3.5}}\smile\put(-3.4,3.5){\line(1,0){3.5}}}
\big (  \grave{\varepsilon}(\overline{\quad\; }\frak{g}(\varepsilon))= 
{\boldsymbol\forall} 
  \rightarrow \frak{g}({\boldsymbol\forall}) \big )  \Bigg ),
\end{multline}  
which is obviously true, because of the subcomponent ``$ 
\put(0,2.3){\line(1,0){8}}\put(4,0){\line(0,1){2.3}}\quad 
({\boldsymbol\forall}
=
{\boldsymbol\forall} )$"  and on that very account can never lead to a 
contradiction. } (\cite{gga}, II 262-263, for ($\beta$ and ($\epsilon$ see 
ibid. 256) 
\end{quote}
We  show now that  by basic laws (V) and (III) we can 
achieve immediately
$$  \put(0,2.3){\line(1,0){8}}\put(4,0){\line(0,1){2.3}}\quad 
({\boldsymbol\forall}
= {\boldsymbol\forall} )$$ 
and  
$$\put(0,2.3){\line(1,0){8}}\put(4,0){\line(0,1){2.3}}\quad \bigg ( 
{\boldsymbol\forall}
 = \grave{\varepsilon}\Big ( 
\put(-3.4,3.5){\line(1,0){3.4}}\put(0.1,0.9){\line(0,4){2.3}}\overset{\frak{g}}{\put(-0.2,3.5){\line(1,0){3.5}}\smile\put(-3.4,3.5){\line(1,0){3.5}}}
\big ( \, \grave{\varepsilon}(\overline{\quad\; }\frak{g}(\varepsilon))= 
\varepsilon \,  \rightarrow \frak{g}(\varepsilon)\big ) \Big ) \bigg ) 
$$  
without applying V'c.

By the definition of the identity sign (\cite{gga}, I 11), the 
notation for the course of value (\cite{gga}, I 15) and  $*)$, the symbols 
``${\boldsymbol\forall}$"
and ``$\grave{\varepsilon}\Big (\put(-3.4,3.5){\line(1,0){3.4}}\put(0.1,0.9){\line(0,4){2.3}}\overset{\frak{g}}
{\put(-0.2,3.5){\line(1,0){3.5}}\smile\put(-3.4,3.5){\line(1,0){3.5}}}
\big ( \, \grave{\varepsilon}(\overline{\quad\; }\frak{g}(\varepsilon))= 
\varepsilon \,  \rightarrow \frak{g}(\varepsilon)\big ) \Big ) $"
denote the same object, and, if    
from basic law 
\begin{equation}\tag{V}\label{V}
\vdash (\grave{\varepsilon}g (\varepsilon) = \grave{\alpha}f(\alpha)) = 
( \overset{\frak{a}}{\put(-0.2,3.5){\line(1,0){3.5}}\smile
\put(-3.4,3.5){\line(1,0){3.5}}} \;   g(\frak{a}) = f(\frak{a}))
\end{equation}
 we draw
\begin{multline}\label{v01}
\Bigg ( \grave{\varepsilon}\Big ( 
\put(-3.4,3.5){\line(1,0){3.4}}\put(0.1,0.9){\line(0,4){2.3}}\overset{
\frak{g}}{\put(-0.2,3.5){\line(1,0){3.5}}\smile\put(-3.4,3.5){\line(1,0){3.5}}}
\big ( \, \grave{\varepsilon}(\overline{\quad\; }\frak{g}(\varepsilon)) = 
\varepsilon   \rightarrow \frak{g}(\varepsilon)
\big ) \Big )
= \grave{\varepsilon}( \overline{\quad\; }f(\varepsilon)) \Bigg )  
=  \\
\Bigg ( \text{  } 
\put(-3.4,3.5){\line(1,0){3.4}}\put(0.1,0.9){\line(0,4){2.3}}\overset{\frak{g}}
{\put(-0.2,3.5){\line(1,0){3.5}}\smile\put(-3.4,3.5){\line(1,0){3.5}}}
\big (  \grave{\varepsilon}(\overline{\quad\; }\frak{g}(\varepsilon)) = 
{\boldsymbol\forall}
   \rightarrow \frak{g}({\boldsymbol\forall})  \big ) =
( \overline{\quad\; }f({\boldsymbol\forall})) \Bigg ), 
\end{multline}
which means $ \text{  }
\put(-3.4,3.5){\line(1,0){3.4}}\put(0.1,0.9){\line(0,4){2.3}}\overset{\frak{g}}
{\put(-0.2,3.5){\line(1,0){3.5}}\smile\put(-3.4,3.5){\line(1,0){3.5}}}
\big ( \, \grave{\varepsilon}(\overline{\quad\; }\frak{g}(\varepsilon)) = 
{\boldsymbol\forall}
 \,  \rightarrow \frak{g}({\boldsymbol\forall})  \big ) =
( \overline{\quad\; }f({\boldsymbol\forall})) $ 
is true if and only if 
``$\grave{\varepsilon}\Big ( 
\put(-3.4,3.5){\line(1,0){3.4}}\put(0.1,0.9){\line(0,4){2.3}}\overset{
\frak{g}}{\put(-0.2,3.5){\line(1,0){3.5}}\smile\put(-3.4,3.5){\line(1,0){3.5}}}
\big ( \, \grave{\varepsilon}(\overline{\quad\; }\frak{g}(\varepsilon)) = 
\varepsilon \,  \rightarrow \frak{g}(\varepsilon)
\big ) \Big )$" 
and 
``$\grave{\varepsilon}( \overline{\quad\; }f(\varepsilon))$" 
refer to the same object,
then consequently  
\begin{multline}\label{v1}
 \Bigg ( \grave{\varepsilon}\Big ( 
\put(-3.4,3.5){\line(1,0){3.4}}\put(0.1,0.9){\line(0,4){2.3}}\overset{
\frak{g}}{\put(-0.2,3.5){\line(1,0){3.5}}\smile\put(-3.4,3.5){\line(1,0){3.5}}}
\big ( \, \grave{\varepsilon}(\overline{\quad\; }\frak{g}(\varepsilon)) = 
\varepsilon \,  \rightarrow \frak{g}(\varepsilon)
\big ) \Big )
= {\boldsymbol\forall} \Bigg )  =  \\
\Bigg ( \text{  }
\put(-3.4,3.5){\line(1,0){3.4}}\put(0.1,0.9){\line(0,4){2.3}}\overset{\frak{g}}
{\put(-0.2,3.5){\line(1,0){3.5}}\smile\put(-3.4,3.5){\line(1,0){3.5}}}
\big ( \, \grave{\varepsilon}(\overline{\quad\; }\frak{g}(\varepsilon)) = 
{\boldsymbol\forall}
 \,  \rightarrow \frak{g}({\boldsymbol\forall})  \big ) =
( \overline{\quad\; }f({\boldsymbol\forall}))  \Bigg ).
\end{multline}
 Now substituting
\begin{align}
\tag*{$vi)$}
``\overset{\frak{g}}{\put(-0.2,3.5){\line(1,0){3.5}}\smile\put(-3.4,3.5){\line(1,0){3.5}}}
\big ( \,\grave{\varepsilon}(\overline{\quad\; }\frak{g}(\varepsilon))= \xi
 \,  \rightarrow \frak{g}(\xi)\big )" 
&\text{ for } ``f(\xi)"
\end{align}
in ($7$) we obtain
\begin{multline}\label{v2}
\Bigg ( \grave{\varepsilon}\Big ( 
\put(-3.4,3.5){\line(1,0){3.4}}\put(0.1,0.9){\line(0,4){2.3}}\overset{
\frak{g}}{\put(-0.2,3.5){\line(1,0){3.5}}\smile\put(-3.4,3.5){\line(1,0){3.5}}}
\big (  \grave{\varepsilon}(\overline{\quad\; }\frak{g}(\varepsilon))= 
\varepsilon   \rightarrow \frak{g}(\varepsilon)
\big ) 
\Big ) 
= {\boldsymbol\forall} \Bigg )  = \\
\Bigg ( \text{  }
\put(-3.4,3.5){\line(1,0){3.4}}\put(0.1,0.9){\line(0,4){2.3}}\overset{\frak{g}}
{\put(-0.2,3.5){\line(1,0){3.5}}\smile\put(-3.4,3.5){\line(1,0){3.5}}}
\big ( \, \grave{\varepsilon}(\overline{\quad\; }\frak{g}(\varepsilon)) = 
{\boldsymbol\forall}
 \,  \rightarrow \frak{g}({\boldsymbol\forall}) =  \\
\overset{\frak{g}}{\put(-0.2,3.5){\line(1,0){3.5}}\smile\put(-3.4,3.5){\line(1,0){3.5}}}
\big ( \grave{\varepsilon}(\overline{\quad\; }\frak{g}(\varepsilon))= 
{\boldsymbol\forall}
   \rightarrow \frak{g}({\boldsymbol\forall})\big ) \Bigg ).
\end{multline}
We replace then
\begin{align}
\tag*{$vii)$}
``\overset{\frak{g}}{\put(-0.2,3.5){\line(1,0){3.5}}\smile\put(-3.4,3.5){\line(1,0){3.5}}}
\big ( \,\grave{\varepsilon}(\overline{\quad\; }\frak{g}(\varepsilon))= \xi
 \,  \rightarrow \frak{g}(\xi)\big )" 
&\text{ for } ``g(\xi)",\\
\tag*{$viii)$}
 ``{\boldsymbol\forall}" & \text{ for } ``a", \\
\tag*{$ix)$}
 ``{\boldsymbol\forall}" & \text{ for } ``b"
\end{align}
in ($2$) yielding
\begin{multline}\label{xb}
\Bigg ( 
\put(-3.4,3.5){\line(1,0){3.4}}\put(0.1,0.9){\line(0,4){2.3}}\overset{\frak{g}}{\put(-0.2,3.5){\line(1,0){3.5}}\smile\put(-3.4,3.5){\line(1,0){3.5}}}
\big ( \,
 \grave{\varepsilon}(\overline{\quad\; }\frak{g}(\varepsilon))= 
{\boldsymbol\forall}
 \,  \rightarrow \frak{g}({\boldsymbol\forall})
\big )  \: = \\
\;\overset{\frak{g}}{\put(-0.2,3.5){\line(1,0){3.5}}\smile\put(-3.4,3.5){\line(1,0){3.5}}}
\big ( \,
 \grave{\varepsilon}(\overline{\quad\; }\frak{g}(\varepsilon))= 
{\boldsymbol\forall}
 \,  \rightarrow 
 \frak{g}({\boldsymbol\forall})
\big ) \Bigg )  \rightarrow \\
\put(0,2.3){\line(1,0){8}}\put(4,0){\line(0,1){2.3}}\quad ( 
{\boldsymbol\forall}
= {\boldsymbol\forall} ).
\end{multline} 
Hence for $*)$, ($8$) and ($9$) 
\begin{equation}\label{noidentical}
 \put(0,2.3){\line(1,0){8}}\put(4,0){\line(0,1){2.3}}\quad ( 
{\boldsymbol\forall}
= {\boldsymbol\forall} ).
\end{equation}

Let us regard everything by means of a different substitution. 
By (V) we can 
yield
\begin{multline}\label{vv1}
 \Bigg ( \grave{\varepsilon}\Big ( 
\put(-3.4,3.5){\line(1,0){3.4}}\put(0.1,0.9){\line(0,4){2.3}}\overset{
\frak{g}}{\put(-0.2,3.5){\line(1,0){3.5}}\smile\put(-3.4,3.5){\line(1,0){3.5}}}
\big ( \, \grave{\varepsilon}(\overline{\quad\; }\frak{g}(\varepsilon))= 
\varepsilon \,  \rightarrow \frak{g}(\varepsilon)
\big ) 
\Big ) 
= {\boldsymbol\forall} \Bigg ) =\\
\Bigg ( \text{  }
\put(-3.4,3.5){\line(1,0){3.4}}\put(0.1,0.9){\line(0,4){2.3}}\overset{\frak{g}}
{\put(-0.2,3.5){\line(1,0){3.5}}\smile\put(-3.4,3.5){\line(1,0){3.5}}}
\Bigg (  \grave{\varepsilon}(\overline{\quad\; }\frak{g}(\varepsilon)) = 
\grave{\varepsilon}\Big ( 
\put(-3.4,3.5){\line(1,0){3.4}}\put(0.1,0.9){\line(0,4){2.3}}\overset{
\frak{g}}{\put(-0.2,3.5){\line(1,0){3.5}}\smile\put(-3.4,3.5){\line(1,0){3.5}}}
\big (  \grave{\varepsilon}(\overline{\quad\; }\frak{g}(\varepsilon))= 
\varepsilon   \rightarrow \frak{g}(\varepsilon)
\big ) 
\Big ) 
 \rightarrow \\
 \frak{g}\Big (
\grave{\varepsilon}\Big (
\put(-3.4,3.5){\line(1,0){3.4}}\put(0.1,0.9){\line(0,4){2.3}}\overset{
\frak{g}}{\put(-0.2,3.5){\line(1,0){3.5}}\smile\put(-3.4,3.5){\line(1,0){3.5}}}
\big (  \grave{\varepsilon}(\overline{\quad\; }\frak{g}(\varepsilon))= 
\varepsilon   \rightarrow \frak{g}(\varepsilon)
\big ) 
\Big ) \Big) \Bigg ) 
 =
( \overline{\quad\; }f({\boldsymbol\forall}))
   \Bigg ). 
\end{multline}
Then for $vi)$ in ($11$), and  $vii)$, $ix)$ and
\begin{align}
 `` \grave{\varepsilon}\Big ( 
\put(-3.4,3.5){\line(1,0){3.4}}\put(0.1,0.9){\line(0,4){2.3}}\overset{\frak{g}}{\put(-0.2,3.5){\line(1,0){3.5}}\smile\put(-3.4,3.5){\line(1,0){3.5}}}
\big ( \, \grave{\varepsilon}(\overline{\quad\; }\frak{g}(\varepsilon))= 
\varepsilon \,  \rightarrow \frak{g}(\varepsilon)
\big ) 
\Big )" \text{ for } ``a",\tag*{$x)$}
\end{align}
in ($2$)  we achieve 
\begin{equation}\label{noidentical1} 
 \put(0,2.3){\line(1,0){8}}\put(4,0){\line(0,1){2.3}}\quad \bigg ( 
{\boldsymbol\forall}
 = \grave{\varepsilon}\Big ( 
\put(-3.4,3.5){\line(1,0){3.4}}\put(0.1,0.9){\line(0,4){2.3}}\overset{\frak{g}}{\put(-0.2,3.5){\line(1,0){3.5}}\smile\put(-3.4,3.5){\line(1,0){3.5}}}
\big ( \, \grave{\varepsilon}(\overline{\quad\; }\frak{g}(\varepsilon))= 
\varepsilon \,  \rightarrow \frak{g}(\varepsilon)
\big ) 
\Big ) \bigg ).
\end{equation}
Accordingly, the above mentioned ($\beta$) and ($\epsilon$) can not be 
drawn, without 
applying (V'c)  or (V'b), and no contradiction can be derived.
We have thus proved that that special self-reference 
procedure which is Russell's antinomy  does not  damage Frege's 
system, QED.
\subsection*{Ending Notes}
Following Frege's \emph{way out},  
all  subsequent literature 
states the necessity to assert the comprehension schema 
 in some restricted form  such as
\begin{equation}\notag  
a \in  \grave{\varepsilon}f(\varepsilon) \leftrightarrow ( a \neq  
\grave{\varepsilon}f(\varepsilon) \wedge 
f(a)),\end{equation}
(see for example \cite{sluga}, 168). On the bases of our 
enlightenment  the uniqueness  established by extensionality  turns out to be 
sufficient 
to prevent the inference of Russell's antinomy, preserving in addition the 
assumption of sound self-reference procedures. As just shown  
we attain (\ref{noidentical}) and (\ref{noidentical1}) as 
consequences of Basic Law (\ref{blIII}), which means that the class of 
classes not belonging to themselves is  
defined by a function which can not take as argument its own course of 
value, i.e. it is a class whose classes are not identical to the class itself.

It would be worthless mentioning all the  well-known attempts to amend the 
foundation of  mathematics
 from the supposed damage to the principle of comprehension, like for 
example 
  Russell's theory of types, Quine's  New Foundation, von Neumann-Bernays  
   system of set theory, and so on. We conclude 
 simply observing  that  in our perspective 
 the necessity of a restriction like Zermelo's 
\emph{Aussonderung} (\cite{foundations}, 36)  appears to be doubtful. Why  
define a Axiom  for 'separating' or 'selecting' those 
members  which fulfil condition $\varphi (x)$, in \emph{Comprehension}, if by 
\emph{Extensionality} $\exists ! y \forall x (x \in y \leftrightarrow 
\varphi (x))$?

\end{document}